\documentclass[draft,english]{lematema}

\usepackage{cleveref}
\usepackage{amsmath,amssymb,amsfonts,amsthm}
\usepackage{xcolor}

\title{Maximum Likelihood Estimation \\ for Nets of Conics}
\titlemark{ML Estimation for Nets of Conics}

\newcommand{\mld}{\mathrm{mld}}
\newcommand{\rmld}{\mathrm{rmld}}
\newcommand{\cL}{\mathcal{L}}
\newcommand{\adj}{\textnormal{adj}}
\newcommand{\rank}{\textnormal{rank}}
\newcommand{\trace}{\textnormal{trace}}

\newcommand{\C}{\mathbb{C}}

\newcommand{\Ps}{\mathbb{P}}
\renewcommand{\S}{\mathbb{S}}
\DeclareMathOperator{\GL}{GL}

\keywords{Nets of quadrics, linear spaces of symmetric matrices, ML-degrees, Veronese surface}

\author{Stefan Dye}
\address{%
\email{stefandye@gmail.com}
}

\author{Kathlén Kohn}
\address{%
KTH Royal Institute of Technology\\
\email{kathlen@kth.se}
}

\author{Felix Rydell}
\address{%
KTH Royal Institute of Technology\\
\email{felixry@kth.se}
}

\author{Rainer Sinn}
\address{%
Universität Leipzig\\
\email{rainer.sinn@uni-leipzig.de}
}

\begin{document}

\maketitle

\begin{abstract}
We study the problem of maximum likelihood estimation for $3$-di\-men\-sional linear spaces of $3\times 3$ symmetric matrices from the point of view of algebraic statistics where we view these nets of conics as linear concentration or linear covariance models of Gaussian distributions on $\mathbb{R}^3$. In particular, we study the reciprocal surfaces of nets of conics which are rational surfaces in $\Ps^5$. We show that the reciprocal surfaces are projections from the Veronese surface and determine their intersection with the polar nets. This geometry explains the maximum likelihood degrees of these linear models.  We compute the reciprocal maximum likelihood degrees. This work is based on Wall's classification of nets of conics from 1977.
\end{abstract}

\section{Introduction}
The aim of this paper is to take a new look at classical results on nets of conics from the point of view of maximum likelihood degrees. 
Maximum likelihood estimation is a widespread approach to fit empirical data to a statistical model based on maximizing the likelihood function. 

We are interested in the generic number of complex critical points of this optimization problem, which is known as the maximum likelihood degree of the model~\cite{catanese2006maximum}.
The models we study are $3$-dimensional sets of trivariate Gaussian distributions with mean zero
that are linear in the space of covariance or concentration matrices.

The $2$-dimensional linear covariance and linear concentration models (in any number of variables) have been covered in \cite{fevola2020pencils}. 
Our paper is a complete case study for the next interesting case, the $3$-dimensional models associated to linear spaces of symmetric matrices. To get a full picture, including all degenerate models, we restrict to Gaussian distributions in three variables and rely on classical results \cite{wall1977nets} by Wall, classifying them from the geometric point of view. 

\subsection{Main results}
For every net of conics,
we determine its reciprocal surface as well as its maximum likelihood degree and its reciprocal maximum likelihood degree. A \emph{net of conics} is a $3$-dimensional linear subspace of the $6$-dimensional space $\S^3$ of $3\times 3$ symmetric matrices. 

For a net of conics $\cL\subset \mathbb{S}^3$, the \emph{reciprocal surface} $\Ps \cL^{-1}$ is the Zariski closure of the set $\Ps (\{A^{-1}\mid A\in \cL, \, \det(A)\neq 0 \})\subset \Ps^5$, where $A^{-1}$ denotes the inverse as a $3\times 3$ matrix. This makes sense for nets that contain an invertible matrix. We call these nets \emph{regular}. 

The reciprocal surfaces are clearly rational. In fact, we show that they are all projections of the Veronese surface in $\Ps^5$. 
Moreover, the reciprocal surface of a net of conics only depends on the type of the net up to congruence action and we identify the center of the projection for every type of regular net.

We also determine the maximum likelihood degree $\mld(\mathcal{L})$ of every regular net $\mathcal{L} \subset \S^3$. 
In addition, we compute using \texttt{Macalauy2}~\cite{m2} the reciprocal maximum likelihood degree $\rmld(\mathcal{L})$
which is the maximum likelihood degree of $\mathcal{L}^{-1}$.
Our results are summarized in \Cref{tab:mldegs}.

Our main tool here is the study of the \emph{ML-base locus} which is the intersection of the reciprocal surface $\Ps \cL^{-1}$ with the polar net $\Ps \cL^\perp$ defined via the trace pairing.
The ML-base loci for the various types of regular nets are summarized in \Cref{tab:MLbase}. We discuss the similarities between the columns of the table in \Cref{sec:MLdegree}.

In \cite{fevola2020pencils}, the authors observed that the intriguing relation
\begin{align}
\label{eq:rmldRelation}
    \rmld(\cL) = \deg(\Ps\cL^{-1}) + \mld(\cL) - 1
\end{align}
holds for pencils $\cL \subset \S^n$ of quadrics. Our results show that this does not generalize to nets of conics.  In fact, neither inequality holds for all regular nets (see the last row of Table 1). However, it does hold for generic nets of conics.
More generally, by comparing the two tables in \cite[Table 1]{linearCovariance}, we see that~\eqref{eq:rmldRelation} holds for generic pencils, for generic linear spaces of symmetric $2 \times 2$ or $3 \times 3$ matrices, but in general \emph{not} for other cases of generic linear spaces $\cL\subset \S^n$.

\begin{table}[ht]
\centering
\begin{tabular}{c|ccccccccccccc}
Type & $A$ & $B$ & $B^*$ & $C$  & $D$ & $D^*$ & $E$ & $E^*$ & $F$ & $F^*$ & $G$ & $G^*$ & $H$\\ \hline
Codim & $0$ & $1$ & $1$ & $2$ & $2$ & $2$ & $3$ & $3$ & $3$ & $3$ & $4$ & $4$ & $5$\\  \hline
$\deg \Ps \mathcal{L}^{-1}$ & 4 & 3 & 4 & 3 & 2 & 4 & 1 & 4 & 2 & 2 & 1 & 2 & 1\\ 
mld $\mathcal{L}$ & 4 & 3 & 3 & 2 & 2 & 2 & 1 & 1 & 0 & 1 & 0 & 0 & 0\\ 
rmld $\mathcal{L}$ & 7 & 5 & 6 & 4 & 3 & 5 & 1 & 4 & 2 & 1 & 0 & 1 & 0\\\hline
Relation~\eqref{eq:rmldRelation} & = & = & = & = & = & = & = & = & $>$ & $<$ & = & = & =\\
\end{tabular}
\caption{For every type of regular net, we list codimension of the set of all nets of that type in the Grassmanian $\mathrm{Gr}(3,\mathbb{S}^3)$,
 degree of the reciprocal surface, ML-degree, reciprocal ML-degree,
 and inequality between the actual reciprocal ML-degree and the one predicted by~\eqref{eq:rmldRelation}.
 }
\label{tab:mldegs}
\end{table}

\begin{table}
\centering
\begin{tabular}{ c | c | c  }
Type & ML-base locus & Common zeroes of the conics \\ \hline
$A$ & $\emptyset$ & $\emptyset$\\  
$B$ & $\emptyset$ & $\emptyset$\\
$B^*$ & $z^2$ & $\{(0:0:1)\}$\\
$C$ & $z^2$ & $\{(0:0:1)\}$\\
$D$ & $\emptyset$ & $\emptyset$\\
$D^*$ & $x^2,y^2$ & $\{(1:0:0),(0:1:0)\}$\\
$E$ & $\emptyset$ & $\emptyset$\\
$E^*$ & $x^2,y^2,z^2$ & $\{(0:0:1),(0:1:0),(1:0:0)\}$\\
$F$ & $z^2$ (mult. 2) & $\{(0:0:1)\}$ (mult. 2)\\
$F^*$ & $y^2-z^2$ & $\{(0:1:i),(0:1:-i)\}$\\
$G$ & $z^2$ & $\{(0:0:1)\}$\\
$G^*$ & $z^2$ (mult. 2) & $\{(0:1:0),(0:0:1)\}$\\
$H$ & $\mathrm{span}\{yz,z^2\}$ & $\{(0:0:1)\}$\\
\end{tabular}
\caption{The table shows the ML-base locus $\Ps \cL^{-1} \cap \Ps \cL^\perp$ and the common zeroes of the nets $\cL$ for the representatives of the types given in \Cref{tab:wallsclass}.}
\label{tab:MLbase}
\end{table}

\section{Wall's classification of nets}

C.~T.~C.~Wall found that any complex net of conics can be categorized as one of 15 geometric types that Wall refers to as $A$, $B$, $B^*$, $C$, $D$,  $D^*$, $E$, $E^*$, $F$, $F^*$, $G$, $G^*$, $H$, $I$, $I^*$. 
Wall's main tool in distinguishing those types is the cubic discriminant curve  obtained by intersecting the net $\Ps \cL$ with the determinantal hypersurface in $\Ps \S^3$.
Nets of type $A$ are exactly those whose discriminant is a smooth (reduced) cubic curve.
In particular, generically chosen nets are of type $A$.
For the geometric description of the other types we refer to Wall~\cite{wall1977nets}.
Bases for nets of each type are given in Table~\ref{tab:wallsclass}. 

\begin{table}[ht]
\centering
\begin{tabular}{ c| c c c c }
Type & $S_1$ & $S_2$ & $S_3$ \\ \hline
$A$  & $y^2+2xz$ & $2yz$ & $-x^2-2gy^2+cz^2+2gxz$ \\  
$B$  & $y^2+2xz$ &  $2yz$&$-x^2-2y^2-9z^2+2xz$ \\
$B^*$ & $y^2+2xz$ & $2yz$&$-x^2-2y^2+2xz$\\
$C$  & $y^2+2xz$ & $2yz$&$-x^2$ \\
$D$  & $x^2$ & $y^2$ & $z^2+2xy $\\
$D^*$ & $2xz$ & $2yz$& $z^2+2xy $\\
$E$  & $x^2$ & $y^2$& $z^2 $\\
$E^*$ & $2xz$ & $2yz$& $2xy $\\
$F$  & $x^2$ & $y^2$& $2xz +2yz$\\
$F^*$ & $x^2$ & $2xy$& $y^2+z^2 $\\
$G$  & $x^2$ & $y^2$& $2yz $\\
$G^*$ & $x^2$ & $2xy$& $2yz $\\
$H$  & $x^2$ &  $2xy$ & $y^2+2xz$\\
$I$  & $x^2$ &  $2xy$ & $y^2$\\
$I^*$ & $2xz$ &  $2yz$ & $z^2$\\

\end{tabular}
\caption{Wall's types of nets of conics with a choice of generators $S_1,S_2,S_3$. Here, $0\neq c\neq -9g^2$ are arbitrary constants.
}\label{tab:wallsclass}
\end{table}

The group $\GL(n)$ acts on the space $\S^n$ of symmetric matrices by congruence:
$\GL(n) \times$ $\mathbb{S}^{n} \rightarrow \mathbb{S}^{n}$, $(g,M)\rightarrow g^TMg$.
Every type, except the generic type $A$, is one orbit of nets under the congruence action. Nets of type $A$ are an irreducible, one-dimensional family of orbits, see \cite[p.~359]{wall1977nets}. This fact follows almost immediately from Wall's construction.

The most degenerate types of nets $I$ and $I^*$ are \emph{singular}, i.e. they do not contain any matrix of full rank. We therefore do not consider them further in this article.

On the ${\binom{n+1}{2}}$-dimensional vector space of complex symmetric $n\times n$ matrices we fix the bilinear trace pairing $(M,N)\mapsto M\bullet N = \trace(MN)$. Restricted to the real symmetric matrices, this is an inner product. For any linear subspace $\cL\subset \S^n$, we write \[
\cL^\perp = \{N\in \S^n\mid M\bullet N = 0 \text{ for all } M\in \cL\}
\]
for its annihilator with respect to the chosen pairing.
The following lemma shows that equivalence classes of nets $\cL \subset \S^3$ under the congruence action
correspond to equivalence classes of the polar nets $\cL^\perp$.
More specifically, the polar net of a net of type $B$ is of type $B^*$.
Similarly, the types $D$ and $D^*$, $E$ and $E^*$, $F$ and $F^*$, as well as $G$ and $G^*$ are polar to each other.
We note that the types $A$, $C$ and $H$ are self-polar in that sense. 

\begin{lemma} 
\label{lem:congruentPolars}
Let $\cL_1$ and $\cL_2$ be $r$-dimensional vector spaces of symmetric matrices of size $n\times n$. Then $\cL_1$ and $\cL_2$ are equivalent under the congruence action if and only if $\cL_1^{\perp}$ and $\cL_2^{\perp}$ are equivalent as well. 
\end{lemma}

\begin{proof} 
For an invertible matrix $g$, assume that $\cL_1^{\perp}=g^T\cL_2^{\perp}g$. We get $$\cL_1=(\cL_1^{\perp})^{\perp}=(g^T\cL_2^{\perp}g)^{\perp}=\{M\in \mathbb{S}^n \mid M\bullet (g^TNg) = 0 \textnormal{ for all }N\in \cL_2^\perp\}.$$
Since $M\bullet (g^TNg)=\trace(Mg^TNg)=\trace(gMg^TN)=(gMg^T)\bullet N$, 
we see that
$\cL_1 = \{ M \in \S^n \mid gMg^T \in \cL_2 \} =g^{-1}\cL_2g^{-T}$, which shows one direction. The other direction is proven analogously.
\end{proof}


\section{Reciprocal surfaces}
\label{sec:recipSurfaces}

We are interested in the study of the reciprocal surface $\Ps \cL^{-1} \subset \Ps^5$ of an arbitrary regular net of conics $\cL$. 
We first observe that it is sufficient to study the reciprocal surface of one representative per congruence class of nets.

\begin{prop}
\label{reciprocalcongruence}
Let $\cL_1, \cL_2 \subset \S^n$ be regular linear subspaces.
If $\cL_1$ and $\cL_2$ are equivalent under the congruence action, then their reciprocal varieties $\Ps \cL_1^{-1}$ and $\Ps \cL_2^{-1}$ in $\Ps \S^n$  are projectively equivalent, even congruent. 
\end{prop}

\begin{proof}
Suppose $\cL_2 = g^T \cL_1 g$ where $g \in \mathrm{GL}(n)$ is an invertible matrix. 
We consider the Zariski open and dense subsets $U_i := \Ps \{ M^{-1} \mid M \in \cL_i, \rank(M)=n \}$ of $\Ps \cL_i^{-1}$ for $i = 1,2$.
We see that $U_2 = g^{-1} U_1 g^{-T}$.
In other words, 
the automorphism 
$\varphi: N \mapsto g^{-1} N g^{-T}$
on $\Ps \S^n$ maps $U_1$ to $U_2$, i.e.
$\varphi(U_1)=U_2$.
Since the map $\varphi$ is closed~\cite[Section 5.8, Theorem 6]{cox2013ideals}, we have that $\varphi(\Ps \cL_1^{-1}) = \Ps \cL_2^{-1}$.
\end{proof}

The reciprocal surfaces of regular nets of conics are all closely related to the Veronese surface in $\Ps^5$ which is the image of the embedding
\begin{align*}
    \nu: \quad\quad\quad \Ps^2 &\longrightarrow \Ps^5, \\
(\alpha:\beta:\gamma)&\longmapsto 
\begin{bmatrix}
\alpha^2 & \alpha\beta & \alpha\gamma \\
\alpha\beta & \beta^2 & \beta\gamma \\
\alpha\gamma & \beta\gamma & \gamma^2
\end{bmatrix}\cong (\alpha^2:\beta^2:\gamma^2:\alpha\beta:\alpha\gamma:\beta\gamma).
\end{align*}

\begin{thm}\label{thm:reciprocalsurfaceA}
For every net $\cL$ of types $A$, $B^*$, $D^*$ and $E^*$, the reciprocal surface $\Ps\cL^{-1}$ is projectively equivalent to the Veronese surface in $\mathbb{P}^5$.
\end{thm}

\begin{proof}
We prove the theorem for types $A$ and $B^*$. The cases of $D^*$ and $E^*$ are performed similarly. 
According to Wall, every net of type $A$ or $B^*$ is congruent to a net spanned by $S_1,S_2,S_3$ as in the first two rows of Table~\ref{tab:wallsclass}.
An arbitrary element of a net with that basis is
$$\alpha S_1+\beta S_2+\gamma S_3=
\begin{bmatrix}
x & y & z
\end{bmatrix}\begin{bmatrix}-\gamma& 0 & \alpha+g\gamma\\
0& \alpha-2g\gamma & \beta \\
\alpha+g\gamma& \beta & c\gamma\end{bmatrix}\begin{bmatrix}
x \\ y \\ z
\end{bmatrix},$$
where $c=0$ and $g=1$ for type $B^*$, and $0 \neq c \neq -9g^2$ for type $A$.
 We consider the adjugate matrix
$$\begin{bmatrix}c\gamma(\alpha-2g\gamma)-\beta^2& \beta(\alpha+g\gamma) & -(\alpha+g\gamma)(\alpha-2g\gamma)\\
\beta(\alpha+g\gamma)& -c\gamma^2-(\alpha+g\gamma)^2 & \beta\gamma \\
-(\alpha+g\gamma)(\alpha-2g\gamma)& \beta\gamma & -\gamma(\alpha-2g\gamma)\end{bmatrix}. $$
The reciprocal surface is the (Zariski closure of the) image of this adjugate map because $\det(M) M^{-1} =  \adj(M)$ for invertible matrices. 
We express the adjugate map in terms of the quadratic monomials in the variables $\alpha,\beta,\gamma$ as such
\begin{align}
\label{eq:transformationMatrix}
\begin{bmatrix}
c\gamma(\alpha-2g\gamma)-\beta^2 \\ -c\gamma^2-(\alpha+g\gamma)^2 \\
 -\gamma(\alpha-2g\gamma) \\ \beta(\alpha+g\gamma) \\ -(\alpha+g\gamma)(\alpha-2g\gamma)
 \\ \beta\gamma 
\end{bmatrix}
=
\begin{bmatrix}0 & -1 & -2cg& 0 & c& 0\\ -1 & 0 & -(c+g^2)& 0 & -2g& 0\\
0 & 0 & 2g& 0 & -1& 0\\
0 & 0 & 0& 1 & 0& g\\
-1 & 0 & 2g^2& 0 & g& 0\\
0 & 0 & 0& 0 & 0& 1\end{bmatrix}\begin{bmatrix}\alpha^2\\ \beta^2\\ \gamma^2\\ \alpha\beta\\\alpha\gamma\\ \beta\gamma\end{bmatrix}.
\end{align}
The determinant of the $6 \times 6$ matrix is equal to $9g^2+c$. For types $A$ and $B^*$, this constant is non-zero and therefore the reciprocal surface is projectively equivalent to the Veronese surface in $\Ps^5$.
\end{proof}

We call the $6\times 6$ matrix in~\eqref{eq:transformationMatrix} the \emph{transformation matrix} corresponding to the reciprocal surface.
For an arbitrary net $\cL$, it is calculated with the same procedure as above.
For the nets of types $B$ and $C$ in \Cref{tab:wallsclass}, the transformation matrix is the same as in \eqref{eq:transformationMatrix}
with $c=-9$ and $g=1$ for type $B$,
and $c=0=g$ for type $C$, respectively.
The transformation matrices for the remaining nets in \Cref{tab:wallsclass} are as follows:

{\footnotesize
$$\begin{array}{c|c|c|c}
    D & D^* & E & E^* \\ 
    \left[\begin{smallmatrix} 0 & 0 & 0 & 0& 0 & 1 \\ 0&0&0&0&1&0\\0&0&-1&1&0&0\\ 0 & 0 & -1 &0&0&0 \\ 0 & 0 & 0&0&0&0 \\ 0&0&0&0&0&0\end{smallmatrix}\right] & 
    \left[\begin{smallmatrix} 0 & -1 & 0 & 0& 0 & 0 \\ -1&0&0&0&0&0\\0&0&-1&0&0&0\\ 0 & 0 & -1 &1&0&0 \\ 0 & 0 & 0&0&0&1 \\ 0&0&0&0&1&0\end{smallmatrix}\right] &
    \left[\begin{smallmatrix} 0 & 0 & 0 & 0& 0 & 1\\ 0&0&0&0&1&0 \\0&0&0&1&0&0\\ 0 & 0 & 0 &0&0&0 \\ 0 & 0 & 0&0&0&0 \\ 0&0&0&0&0&0\end{smallmatrix}\right] &
    \left[\begin{smallmatrix} 0 & -1 & 0 & 0& 0 & 0 \\ -1&0&0&0&0&0\\0&0&-1&0&0&0\\ 0 & 0 & 0 &1&0&0 \\ 0 & 0 & 0&0&0&1 \\ 0&0&0&0&1&0\end{smallmatrix}\right]
\end{array}$$

$$\!\!\!
\begin{array}{c|c|c|c|c}
   F & F^* & G & G^* & H  \\
   \left[\begin{smallmatrix} 0 & 0 & -1 & 0& 0 & 0 \\ 0&0&-1&0&0&0\\0&0&0&1&0&0\\ 0 & 0 & 1 &0&0&0 \\ 0 & 0 & 0&0&0&-1 \\ 0&0&0&0&-1&0\end{smallmatrix}\right] &
   \left[\begin{smallmatrix} 0 & 0 & 1 & 0& 0 & 0 \\ 0&0&0&0&1&0\\0&-1&0&0&1&0\\ 0 & 0 & 0 &0&0&-1 \\ 0 & 0 & 0&0&0&0 \\ 0&0&0&0&0&0\end{smallmatrix}\right] &
   \left[\begin{smallmatrix} 0 & 0 & -1 & 0& 0 & 0 \\ 0&0&0&0&0&0\\0&0&0&1&0&0\\ 0 & 0 & 0 &0&0&0 \\ 0 & 0 & 0&0&0&0 \\ 0&0&0&0&-1&0\end{smallmatrix}\right] &
   \left[\begin{smallmatrix} 0 & 0 & -1 & 0& 0 & 0 \\ 0&0&0&0&0&0\\0&-1&0&0&0&0\\ 0 & 0 & 0 &0&0&0 \\ 0 & 0 & 0&0&0&1 \\ 0&0&0&0&-1&0\end{smallmatrix}\right] &
   \left[\begin{smallmatrix} 0 & 0 & 0 & 0& 0 & 0\\ 0&0&-1&0&0&0 \\0&-1&0&0&1&0\\ 0 & 0 & 0 &0&0&0 \\ 0 & 0 & -1&0&0&0 \\ 0&0&0&0&0&1\end{smallmatrix}\right]
\end{array}\!\!\!$$
}



An immediate consequence from our calculations is the following fact that will be useful for computing maximum likelihood degrees in the next section.

\begin{cor} 
\label{cor:noRankOne}
Nets of types $A$, $B^*$, $D^*$ and $E^*$ have no matrices of rank 1.
\end{cor}

\begin{proof}
The adjugate of a $3\times 3$ symmetric matrix of rank $1$ is the zero matrix, but the $6\times 6$ transformation matrix corresponding to a net of type $A$, $B^*$, $D^*$, or $E^*$ is invertible, which shows that the adjugate of any non-zero matrix in the net is non-zero.
\end{proof}

\begin{rem}
\label{rem:rankOneZeroLocus}
We could have derived Corollary~\ref{cor:noRankOne} from Wall's classification of nets. From that statement, we can derive Theorem~\ref{thm:reciprocalsurfaceA} by a geometric argument: The adjugate map restricted to nets of type $A$, $B^*$, $D^*$, and $E^*$ gives a morphism from the net to $\Ps^5$ defined globally by quadratic forms. Since the image is non-degenerate, the image is the Veronese surface as claimed. 

More generally, there is a geometric explanation in the background of our computations using polarity: A matrix of rank $1$ in a net $\cL$ is a base point of the polar net $\cL^\perp$ by the usual trace trick \[ M\bullet (vv^T) = \trace(M (vv^T)) = v^T M v. \]
We chose to take the computational road to avoid more careful geometric arguments with multiplicities for the more degenerate nets that we discuss next.
\end{rem}

\begin{prop}
\label{prop:projectionsDegreeThree} The reciprocal surfaces of nets of types $B$ and $C$ have degree $3$ and are projections of the Veronese surface from a  point on it.
\end{prop}

\begin{proof} 
For these types, the transformation matrices do not have full rank. The right-kernel is the center of the projection. A calculation shows that this kernel is $(4g^2:0:1:0:2g:0)$,
which corresponds to the symmetric matrix 
$$\left[\begin{smallmatrix} 4g^2 & 0 & 2g \\ 0 & 0 & 0 \\ 2g & 0 & 1\end{smallmatrix} \right] \cong 4g^2 x^2 + z^2 + 4g xz$$
($g= 1$ for $B$ and $g=0$ for $C$). This point is equal to $\nu(2g:0:1) = (2gx+z)^2$ and therefore lies on the Veronese surface. 
\end{proof}

\begin{prop}
\label{prop:projectionsDegreeTwo}
The reciprocal surfaces of nets of types $D$, $F$, $F^*$, and $G^*$ have degree $2$ and are projections of the Veronese surface from a line. For types $D$ and $F$, the center of projection is a secant line to the Veronese surface (spanned by two distinct points on it).
For types $F^*$ and $G^*$, the line is a tangent line. 
\end{prop}

\begin{proof} Again, we calculate the right-kernel of the transformation matrices that have rank $4$ in these cases.

\noindent
\textbf{Types} $D$ \& $F$:
The kernel is spanned by $\nu(1:0:0)\cong x^2$ and $\nu(0:1:0)\cong y^2$.

\noindent
\textbf{Types} $F^*$ \& $G^*$:
The kernel is spanned by the point $\nu(1:0:0)\cong x^2$, and the point $(0:0:0:1:0:0) \in T_{\nu(1:0:0)} \nu(\Ps^2)$ which corresponds to the tangent direction $2xy$.
\end{proof}

\begin{rem}
\label{rem:singularLocus}
The reciprocal surfaces in Proposition~\ref{prop:projectionsDegreeThree} are rational normal scrolls in $\Ps^4$, which are isomorphic to the Hirzebruch surface $\mathcal{H}_1$. As toric varieties, they correspond to the lattice polygon \[
{\rm conv}\{(0,0),(2,0),(0,1),(1,1)\}.
\]
The reciprocal surfaces in Proposition~\ref{prop:projectionsDegreeTwo} for nets of types $D$ and $F$ are smooth quadrics in $\Ps^3$ and are therefore projectively equivalent to $\Ps^1\times \Ps^1$. As toric varieties, they correspond to the lattice square \[{\rm conv}\{(0,0),(1,0),(0,1),(1,1)\}.\]
For types $F^*$ and $G^*$, the reciprocal surfaces are cones over conics, corresponding to the lattice triangle
\[ {\rm conv}\{(0,0),(2,0),(0,1)\}.\]
In particular, we see that the reciprocal surfaces for types $F^*$ and $G^*$ are the only singular ones.
For the nets $\cL$ in \Cref{tab:wallsclass}, the vertex of the reciprocal cone $\Ps \cL^{-1}$ is $y^2+z^2$ in type $F^*$ and $yz$ in type  $G^*$.
\end{rem}

By inspecting the $6\times 6$ transformation matrices corresponding to the nets of types $E$, $G$ and $H$,
we analogously see that their reciprocal surfaces are planes obtained by projecting the Veronese surface from a plane. 
\begin{itemize}
    \item For type $E$, the projection center is a plane spanned by $\nu(1:0:0) \cong x^2$, $\nu(0:1:0) \cong y^2$, and $(0:0:0:1:0:0) \cong 2xy$, where the latter point corresponds to a tangent direction at the previous two points. That plane intersects the Veronese surface in a conic, so it is a net of type $I$.
    \item For type $G$, the projection center intersects the Veronese surface in two distinct points $\nu(1:0:0)\cong x^2$ and $\nu(0:1:0)\cong y^2$. Moreover, it contains the tangent line in direction $(0:0:0:0:0:1)\cong 2yz$ at the latter of those intersection points. This projection center is a net of type $G$.
    \item For type $H$, the projection center intersects the Veronese surface at a single point $\nu(1:0:0)\cong x^2$. It also contains a tangent at that point in direction $(0:0:0:1:0:0)\cong 2xy$. The remaining generator is $(0:1:0:0:1:0)\cong y^2+2xz$, so it is a net of type $H$.
\end{itemize}

\begin{rem}
\label{rem:annihilators}
For each type, except $A,B^*,D^*$ and $E^*$, the corresponding $6\times 6$ transformation matrix is of rank at most $5$, so it is an embedded projection. 
We compute the left-kernel of each transformation matrix, which precisely defines the hyperplanes that the reciprocal surface lies in. 
In other words, the left-kernel of the transformation matrix is the annihilator of the linear space spanned by the reciprocal surface. The following table provides bases for the annihilators.

\begin{table}[ht]
\centering
\begin{tabular}{ c| c }
Type & Basis for annihilator of $\mathrm{span}\{ \cL^{-1} \}$ \\ \hline  
$B$  & $2xz - y^2 + 3gz^2$ \\
$C$  & $y^2-2xz$ \\
$D$  & $xz$, $yz$ \\
$E$  & $xy,xz,yz $ \\
$F$  & $x^2+2xy$, $y^2+2xy$  \\
$F^*$ &$xz$, $yz$ \\
$G$  & $xy$, $xz$, $y^2$ \\
$G^*$ &  $xy$, $y^2$ \\
$H$  & $x^2$, $xy$, $y^2-2xz $ \\
\end{tabular}
\caption{We list the types whose  $6\times 6$ transformation matrix is singular  and the elements that span its left-kernel (for the generators in \Cref{tab:wallsclass}).
}\label{tab:netannihilators}
\end{table}

For a net $\cL$ of type $E$, we see that its reciprocal surface $\Ps \cL^{-1}$ is a net of type $E$ as well.
Similarly, the reciprocal surface of a net of type $G$ (resp.~$H$) is again a net of type $G$ (resp.~$H$).

\end{rem}


\section{ML-degrees}
\label{sec:MLdegree}
The \emph{maximum likelihood degree (ML-degree)}  of a real linear space $\cL \subset \S^n$ is the number of complex critical points of the log-likelihood function
\begin{align*}
    \ell_S: \cL &\longrightarrow \mathbb{R}, \\
    M &\longmapsto \log\det (M)-\trace (SM) 
\end{align*}
for a generic matrix $S \in \S^n$.
The critical equations are polynomial in the entries of $M$ and therefore the notion of ML-degree also makes sense for complex linear spaces $\cL$.
Our goal is to compute the ML-degree for every net of conics.
The ML-degree is invariant under the congruence action~\cite[Lemma 4.1]{fevola2020pencils}, 
so it is sufficient to determine the ML-degree for every net in \Cref{tab:wallsclass}.
We make use of the following two statements.

\begin{prop}[\cite{mldefective}]
\label{prop:mldefective}
The ML-degree of a linear space $\cL \subset \S^n$ is
at most the degree of its reciprocal variety $\Ps \cL^{-1}$.
This is an equality if and only if the ML-base locus $\Ps \cL^{-1} \cap \Ps \cL^{\perp}$ is empty.
\end{prop}

More specifically, a formula for the ML-degree of a linear space $\cL \subset \S^n$ in terms of Segre classes of its ML-base locus is given in~\cite{mldegLSSMs}.
We will only use the following special case of that formula.

\begin{prop}[\cite{mldegLSSMs}]\label{prop:segreClassFormula} 
For a linear space $\cL \subset \S^n$ whose ML-base locus  is finite and consists only of smooth points of $\Ps \cL^{-1}$, we have
$$\mld (\cL)=\deg( \Ps \cL^{-1})-\deg(\Ps\cL^{-1}\cap \Ps\cL^\perp), $$
where the degree of the ML-base locus is its scheme-theoretic degree (i.e., the constant coefficient of its Hilbert polynomial).
\end{prop}

Hence, the study of the ML-base loci of nets of conics is key in our computation of ML-degrees.
\begin{rem}\label{baselocusfullrank}
We first observe that  ML-base loci do not contain matrices of full rank. Indeed, if a full-rank matrix $M$ would be in the ML-base locus of a linear space $\cL \subset \S^n$, then we have $M^{-1}\in \mathcal{L}$, so $M\in \mathcal{L}^{\perp}$ implies  $0=\trace (M^{-1}M)=\trace (I_n)=n>0$;  a contradiction. 
\end{rem}

\begin{lemma}\label{InvLem} Let $\mathcal{L} \subset \S^n$ be a regular linear space containing only matrices of rank $0$, $n-1$ or $n$. Then $\cL^{-1}$ contains only matrices of rank $0$, $1$ or $n$. 
\end{lemma}

\begin{proof} 
The adjugate map is defined everywhere on $\Ps \cL$, 
so the image of the adjugate map is Zariski closed~\cite[Section 5.8, Theorem 6]{cox2013ideals}. 
Thus, if $M\in \mathcal{L}^{-1}$, then $M=\adj(N)$ for some $N\in \mathcal{L}$. If $N$ is invertible, so is its adjugate. If $N$ is of rank $n-1$, then its adjugate has rank $1$.
\end{proof}

Now we aim to compute the ML-degrees of regular nets.
To make use of Propositions~\ref{prop:mldefective} and \ref{prop:segreClassFormula}, we compute the ML-base locus for every regular net.
By the following lemma, it is sufficient to compute the ML-base locus for one representative per type.

\begin{lemma}\label{congruentMLbase}
Let $\cL_1$ and $\cL_2$ be congruent linear subspaces of $\S^n$. The common zero locus in $\Ps^{n-1}$ of the quadrics in $\cL_1$ is projectively equivalent to the zero locus of $\cL_2$.
Furthermore, the ML-base loci of $\cL_1$ and $\cL_2$ are congruent.
\end{lemma}

\begin{proof} 
Since the congruence action $M \mapsto g^T M g$ for $g\in \GL_n$ corresponds to the change of coordinates $x \mapsto g^{-1}x$ on $\C^n$, the common zero locus of quadrics in $\cL$ and a congruent subspace $g^T\cL g$ are related by a change of coordinates on $\Ps^{n-1}$.

To prove the second part, we relate both $(g^T\cL g)^{-1}$ to $\cL^{-1}$ and $(g^T\cL g)^\perp$ to $\cL^\perp$. So firstly, $(g^T\cL g)^{-1}=g^{-1}\cL^{-1}g^{-T}$ holds, i.e.~the reciprocal varieties of congruent subspaces are congruent, see Proposition~\ref{reciprocalcongruence}. 

Secondly, we also have that $(g^T\cL g)^{\perp}=g^{-1}\cL^{\perp}g^{-T}$, i.e.~the polar subspaces are also congruent, see Lemma~\ref{lem:congruentPolars}.
This proves the claim because this is the same transformation as for the reciprocal varieties. 
\end{proof}

The computation of the ML-base locus for one representative per type is straightforward in \texttt{Macaulay2}. 
The resulting ML-base loci are listed in Table~\ref{tab:MLbase}. The ML-base loci can also be determined by pleasant geometric arguments. We demonstrate how this can be done for the types  $A$, $F$, $F^*$, and $G^*$.

The ML-base locus of a net of conics consists of matrices of rank one or two (see Remark~\ref{baselocusfullrank}).
The rank-one matrices in the ML-base locus correspond to points in $\Ps^2$ in the zero locus of the conics in the net by the trace trick in Remark~\ref{rem:rankOneZeroLocus}.
This explains the apparent similarities between the two columns in \Cref{tab:MLbase}.
However, not every point in the common zero locus of the conics has to yield a point in the ML-base locus.
So to compute the ML-base locus we use only the rank-one matrices in the polar net $\cL^\perp$ that appear in the reciprocal surface, and we find the rank-two matrices in their intersection.

\begin{prop}\label{BigT}
For a net of type $A$, the ML-base locus is empty and the ML-degree is $4$.
\end{prop}

\begin{proof}
It is sufficient to show that the ML-base locus of a net $\cL$ of type $A$ is empty as this implies the assertion by Proposition~\ref{prop:mldefective} and Theorem \ref{thm:reciprocalsurfaceA}. We check that there is no matrix of rank $1,2$ or $3$ in the ML-base locus:\\
\textbf{Rank 1:} Since $\cL^{\perp}$ is also of type $A$, it has no matrix of rank 1 by Corollary~\ref{cor:noRankOne}. 
\\
\textbf{Rank 2:} By Corollary~\ref{cor:noRankOne} and Lemma~\ref{InvLem}, there is no matrix of rank $2$ in $\cL^{-1}$.
\\
\textbf{Rank 3:} By Remark~\ref{baselocusfullrank}, the ML-base locus does not contain full-rank matrices.
\end{proof}

\begin{thm}
The ML-degree of a regular net of conics depends only on its type. All ML-degrees are listed in \Cref{tab:mldegs}. 
\end{thm}

\begin{proof}
The ML-degree is invariant under congruence action by \cite[Lemma 4.1]{fevola2020pencils}.
We have shown that the ML-degree of every net of type $A$, which has infinitely many orbits under congruence, is always $4$ (see Proposition~\ref{BigT}). All other types come in one congruence class, as observed by Wall \cite[p.~359]{wall1977nets}.

To compute the ML-degree for all remaining types besides $A$, we determine the ML-base locus and apply  Propositions~\ref{prop:mldefective} and \ref{prop:segreClassFormula}. We outline here the general strategy and then discuss the types $F$, $F^*$, and $G^*$ in detail below.

After computing the ML-base loci for all types (see \Cref{tab:MLbase}), we observe that $H$ is the only type whose ML-base locus is not finite. By Remark~\ref{rem:annihilators}, the reciprocal surface of a net of type $H$ is a plane, so it has degree one.
Hence, Proposition~\ref{prop:mldefective} shows that the ML-degree for type $H$ is zero.

For all other types besides $H$ the ML-base locus is finite. Moreover, using Remark~\ref{rem:singularLocus} we see that every point in the ML-base locus is a smooth point of the reciprocal surface. 
This allows us to apply Proposition~\ref{prop:segreClassFormula} (together with our results on the degree of reciprocal surfaces in Section~\ref{sec:recipSurfaces}) to compute the ML-degree.
\end{proof}

A helpful tool for performing detailed analysis on nets of different types is $$X(M) := \left\lbrace N \in \S^3 \mid MN=0 \right\rbrace, \quad  \text{ where } M \in \S^n.$$

\begin{lemma}
\label{lem:productZero}
Let $\cL$ be a regular net of conics.
For every $N \in \cL^{-1}$ of rank two, there is a rank-one matrix $M \in \cL$ such that $N \in X(M)$.
\end{lemma}

\begin{proof}
We consider  the Zariski closure of the graph of the matrix inversion map:
\begin{align*}
    \Gamma := \overline{\left\lbrace(M, M^{-1}) \mid M \in \Ps \S^3, \; \rank(M) = 3 \right\rbrace} \subset \Ps \S^3 \times \Ps \S^3.
\end{align*}
It was shown in~\cite{graphMatrixInversion} that
\begin{align}
\label{eq:graph}
    \Gamma = \{ (M,N) \in \Ps \S^3 \times \Ps \S^3 \mid MN = t \cdot I_3 \text{ for some constant } t \}.
\end{align}
Similarly, we consider the Zariski closure $\Gamma_\cL$ of the graph of matrix inversion restricted to the net $\cL$, i.e.
$\Gamma_\cL = \overline{\{(M, M^{-1}) \mid M \in \Ps \cL, \; \rank(M) = 3 \} }$. 
We clearly have the containment $\Gamma_\cL \subset \Gamma$.
The projection $\pi_2: \Gamma_\cL \to \Ps \S^3, (M,N) \mapsto N$ onto the second factor is a morphism whose image is $\Ps \cL^{-1}$.
Hence, for every $N \in \Ps \cL^{-1}$, there is an $M \in \Ps \cL$ with $(M,N) \in \Gamma_\cL$.
In particular, we see from~\eqref{eq:graph} that $MN = t \cdot I_3$ for some constant $t$.
If $N$ is of rank two, then $t$ must be zero and $M$ must be of rank one. 
\end{proof}

\begin{prop}
\label{prop:typeF}
For a net of type $F$, the ML-base locus is a double point of rank $1$ and the ML-degree is $0$.
\end{prop}
\begin{proof}
By Lemma~\ref{congruentMLbase}, it is sufficient to consider the net $\cL$ of type $F$ in \Cref{tab:wallsclass}.
We immediately see that the conics in that net have exactly one point in their common zero locus, namely $(0:0:1)$.
Due to Remark \ref{rem:rankOneZeroLocus}, there is exactly one rank-one matrix in $\Ps \mathcal{L}^{\perp}$, which is
 $N_1 = \left[\begin{smallmatrix} 0 & 0 & 0\\ 0 & 0 & 0\\ 0 & 0 & 1\\ \end{smallmatrix} \right]$.
We find that the adjugate of the rank-two matrix $M_1 = \left[\begin{smallmatrix} 1 & 0 & 0\\ 0 & 1 & 0\\ 0 & 0 & 0\\ \end{smallmatrix} \right] \in \mathcal{L}$  is $N_1$.
This shows that $N_1$ is contained in the ML-base locus.

We will show now that $N_1$ is the only matrix in the ML-base locus.
To prove this we need to exclude the existence of rank-two matrices in the ML-base locus.
We assume for contradiction that a rank-two matrix $N_2$ exists in the ML-base locus. 
By Lemma~\ref{lem:productZero}, 
there exists a rank-one matrix $M_2 \in \mathcal{L}$ such that $N_2 \in X(M_2) \cap \mathrm{span}\{ \cL^{-1} \} \cap \cL^\perp$. 
There are two rank-one matrices in $\Ps \cL$;   $M_{2,1} = \left[\begin{smallmatrix} 1 & 0 & 0\\ 0 & 0 & 0\\ 0 & 0 & 0\\ \end{smallmatrix} \right]$ and 
 $M_{2,2} = \left[\begin{smallmatrix} 0 & 0 & 0\\ 0 & 1 & 0\\ 0 & 0 & 0\\ \end{smallmatrix} \right]$.
However, using Table~\ref{tab:netannihilators}, we see that the intersection $X(M_{2,i}) \cap \mathrm{span}\{ \cL^{-1} \} \cap \cL^\perp$ for $i \in \{ 1,2 \}$ contains only the rank-one matrix $N_1$ (up to scaling). This is a contradiction.
 
So far we have shown that  $N_1$ is the only matrix in the ML-base locus.
Since the surface $\Ps \cL^{-1}$ and the plane $\Ps \cL^\perp$ are both contained in the hyperplane with annihilator $x^2 - y^2$ (see again Table~\ref{tab:netannihilators}), 
they must intersect in $\deg(\Ps \cL^{-1}) = 2$ many points, counted with multiplicity.
Since they intersect in the single point $N_1$, the intersection multiplicity at $N_1$ must be $2$.

Finally, since $\Ps \cL^{-1}$ is smooth, we can apply the formula in Proposition~\ref{prop:segreClassFormula}: 
$\mathrm{mld}(\cL) = 2 - 2 = 0$.
\end{proof}

\begin{prop}
\label{prop:typeFstar}
For a net of type $F^*$, the ML-base locus is a reduced point of rank $2$ and the ML-degree is $1$.
\end{prop}

\begin{proof} Again, it is enough to consider one net $\cL$ of type $F^*$, we choose the one given in \Cref{tab:wallsclass}. With respect to the given generators $S_1,S_2$ and $S_3$, the adjugate map is $$\left[\begin{smallmatrix} \alpha & \beta & 0\\ \beta & \gamma & 0\\ 0 & 0 & \gamma\\ \end{smallmatrix} \right]\longmapsto \left[\begin{smallmatrix} \gamma^2 & -\beta\gamma & 0\\ -\beta\gamma & \alpha\gamma & 0\\ 0 & 0 & \alpha\gamma-\beta^2\\ \end{smallmatrix} \right].$$ 
From this expression, we see directly that the prime ideal of the reciprocal surface  is $ \mathrm{I}(\Ps \cL^{-1}) = \langle u_{13},u_{23},u_{11}(u_{33}-u_{22})+u_{12}^2 \rangle $, where $u_{ij}$ are coordinates for a general symmetric $3\times 3$ matrix. Intersecting the reciprocal surface and the polar net (given by the ideal $\mathrm{I}(\mathbb{P}\cL^{\perp})=\langle u_{11},u_{12},u_{22}+u_{33} \rangle$) yields the ideal $\langle u_{11},u_{12},u_{13},u_{23},u_{22}+u_{33}\rangle$. This ideal is of degree $1$ and  corresponds to a single rank-two matrix. Since $\Ps \cL^{-1}$ is smooth at this point (see  Remark~\ref{rem:singularLocus}), we can apply the formula in Proposition~\ref{prop:segreClassFormula}: 
$\mathrm{mld}(\cL) = 2 - 1 = 1$.
\end{proof} 

\begin{prop}
\label{prop:typeGstar}
For a net of type $G^*$, the ML-base locus is a double point of rank $1$ and the ML-degree is $0$.
\end{prop}

\begin{proof} 
As above, we only need to consider the net $\cL$ of type $G^*$ in \Cref{tab:wallsclass}.
First of all, we have
$$\mathbb{P}\cL^{\perp}=\mathrm{span}\Big\{\left[\begin{smallmatrix} 0 & 0 & 0\\ 0 & 1 & 0\\ 0 & 0 & 0\\ \end{smallmatrix} \right],\left[\begin{smallmatrix} 0 & 0 & 0\\ 0 & 0 & 0\\ 0 & 0 & 1\\ \end{smallmatrix} \right],\left[\begin{smallmatrix} 0 & 0 & 1\\ 0 & 0 & 0\\ 1 & 0 & 0\\ \end{smallmatrix} \right]\Big\}.$$
By Remark~\ref{rem:annihilators}, 
 we  observe:
\begin{align}
\label{eq:baseLocusGstar}
    \mathbb{P}\cL^{\perp}\cap \mathbb{P}\cL^{-1}\subseteq \mathbb{P}\cL^{\perp}\cap \mathrm{span}\{\mathbb{P}\cL^{-1} \} =  \mathrm{span} \Big\{\left[\begin{smallmatrix} 0 & 0 & 0\\ 0 & 0 & 0\\ 0 & 0 & 1\\ \end{smallmatrix} \right],\left[\begin{smallmatrix} 0 & 0 & 1\\ 0 & 0 & 0\\ 1 & 0 & 0\\ \end{smallmatrix} \right]\Big\}.
\end{align}
The only rank-one matrix in the right-hand side is $N_1=\left[\begin{smallmatrix} 0 & 0 & 0\\ 0 & 0 & 0\\ 0 & 0 & 1\\ \end{smallmatrix} \right]$. This matrix lies in the reciprocal surface, since it is the adjugate of $S_2$ (see Table~\ref{tab:wallsclass}). 
We claim that $N_1$ is the only matrix in the ML-base locus. By contradiction, we assume that $N_2$ is a rank-two matrix in the locus. 
Lemma~\ref{lem:productZero}  says that there is a matrix $M\in \cL$ such that $MN_2=0$, and this $M$ must be of rank 1.  The only rank-one matrix in $\Ps \cL$ is $M=\left[\begin{smallmatrix} 1 & 0 & 0\\ 0 & 0 & 0\\ 0 & 0 & 0\\ \end{smallmatrix} \right]$. However, intersecting $X(M)$ with~\eqref{eq:baseLocusGstar} yields only $N_1$, so this contains no matrix of rank two.

Both $\mathbb{P}\cL^{\perp}$ and $\mathbb{P}\cL^{-1}$ are contained in the hyperplane defined by $xy$. Hence, the degree of their intersection is $\deg(\mathbb{P}\cL^{-1})=2$. Since the reciprocal surface $\Ps \cL^{-1}$ is smooth at $N_1$ (see  Remark~\ref{rem:singularLocus}), we can apply the formula in Proposition~\ref{prop:segreClassFormula}: 
$\mathrm{mld}(\cL) = 2 - 2 = 0$.
\end{proof}

For all remaining types the same techniques as used in Propositions~\ref{BigT}, \ref{prop:typeF}, \ref{prop:typeFstar}, and \ref{prop:typeGstar} can be applied to determine the ML-base locus and the ML-degree.

\section{Reciprocal ML-degrees}
The \emph{reciprocal maximum likelihood degree} of a real linear space $\cL \subset \S^n$ is the number of complex critical points of the log-likelihood function
 \begin{align*}
    \ell_S: \cL^{-1} &\longrightarrow \mathbb{R}, \\
    M &\longmapsto \log\det (M)-\trace (SM) 
\end{align*}
for a generic matrix $S \in \S^n$.
Note that in this setting the log-likelihood function is defined on the reciprocal variety $\cL^{-1} \subset \S^n$, instead of on the linear space $\cL$ as in Section~\ref{sec:MLdegree}.
The reciprocal ML-degree is invariant under the congruence action~\cite[Lemma 4.1]{fevola2020pencils}.
So as before, we can determine the reciprocal ML-degree for every type except $A$ by computation with \texttt{Macaulay2}. This is because the nets of every such type form a single orbit under congruence. 
Our results are listed in \Cref{tab:mldegs}.

In type $A$, we have sampled several orbits and always found $7$ using \linebreak[4] \texttt{Macaulay2}. We therefore believe that the reciprocal ML-degree of every net of type $A$ is $7$. 
This agrees with the ``$7$'' found on the left of Table~1 in \cite{linearCovariance} which is the reciprocal ML-degree of generic nets of conics. So we conjecture that the generic nets of conics (generic with respect to the reciprocal ML-degree) are exactly the nets of type $A$ (which are by definition generic with respect to their intersection with the discriminant hypersurface).
We note that the fact that the reciprocal surface of a net of type $A$ is projectively equivalent to the Veronese surface is not sufficient to imply that conjecture because the reciprocal surfaces of nets of types $B^*$, $D^*$, and $E^*$ are also Veronese surfaces (see Theorem~\ref{thm:reciprocalsurfaceA}) but their reciprocal ML-degrees are $6$, $5$, and $4$, respectively.

Further research is required to understand reciprocal ML-degrees using the underlying geometry of linear spaces in $\S^n$ and their reciprocal varieties.

\section*{Acknowledgements}
We thank the organizers of the working group on Linear Spaces of Symmetric Matrices at MPI MiS Leipzig and Tim Seynnaeve for pointing out a useful reference.
KK and FR were supported by the Knut and Alice Wallenberg Foundation within their WASP (Wallenberg AI, Autonomous Systems and Software Program) AI/Math initiative.

\end{document}